\documentclass[letterpaper, 10 pt, conference]{ieeeconf}
\IEEEoverridecommandlockouts                              % This command is only
\usepackage{graphicx}
\usepackage{amsmath}
\usepackage{amssymb}
\usepackage{amsfonts}
\usepackage{latexsym}
\usepackage{epsfig}
\usepackage{theorem}
\usepackage{fancybox}
\usepackage{shadow}
\usepackage{float}
\usepackage{comment}
\usepackage{color}
\usepackage{enumerate}
\usepackage{psfrag}

\newtheorem{propo}{Proposition}
\newtheorem{lem}{\bf{Lemma}}
\newtheorem{theo}{Theorem}

\begin{document}
\title{\LARGE\bf Velocity-aided Attitude Estimation for Accelerated Rigid Bodies}

\author{Minh-Duc Hua, Philippe Martin, Tarek Hamel
\thanks{Minh-Duc Hua is with ISIR UPMC-CNRS (Institute for Intelligent Systems and Robotics), Paris, France, {\tt\footnotesize $hua@isir.upmc.fr$}}
\thanks{Philippe Martin is with Mines ParisTech, PSL Research University, Paris, France,
{\tt\footnotesize $philippe.martin@mines\!\!-\!\!paristech.fr$}}%
\thanks{Tarek Hamel is with I3S UNS-CNRS, Sophia Antipolis, France, $thamel@i3s.unice.fr$}%
\thanks{This work was supported by the French Agence Nationale de la Recherche through the ANR ASTRID SCAR project ``Sensory Control of Aerial Robots'' (ANR-12-ASTR-0033)}
}

\maketitle

\pagestyle{empty}
\thispagestyle{empty}

\begin{abstract}
Two nonlinear observers for velocity-aided attitude estimation, relying on gyrometers, accelerometers, magnetometers, and velocity measured in the body-fixed frame, are proposed. As opposed to state-of-the-art body-fixed velocity-aided attitude observers endowed with local properties, both observers are (almost) globally asymptotically stable, with very simple and flexible tuning. Moreover, the roll and pitch estimates are globally decoupled from magnetometer measurements.
\end{abstract}

\section{Introduction}

Robotic vehicles commonly need to know their orientation and velocity to be operated. When cost or weight is an issue, using very accurate inertial sensors for ``true'' (i.e. based on the Schuler effect due to a non-flat rotating Earth) inertial navigation is excluded. Instead, low-cost systems –-sometimes called velocity-aided Attitude Heading Reference Systems–- rely on light and cheap strapdown gyrometers, accelerometers and magnetometers ``aided'' by velocity sensors (provided for example in body-fixed coordinates by an air-data, a Doppler radar system or a Doppler velocity log (DVL), or in Earth-fixed coordinates by a GPS engine). The various measurements are then ``merged'' according to the motion equations of the vehicle assuming a flat non-rotating Earth, usually with a linear complementary filter or an Extended Kalman Filter (EKF).

Nonlinear attitude observers have become an alternative to the EKF, starting with~\cite{Salcudean91} and then over the last fifteen years~\cite{nf99,ts03,tmrm07,mhp08,vcso08,ms08IFAC,BonnabelITAC08,hua10cep,ms10cep,huaIFAC11,rt11CDC,jensen11,trumpf12,grip12,zamani13,bras2013global,troniICRA13,Dukan13,hua14}. The performance of recent observers is comparable to modern nonlinear filtering techniques~\cite{Crassidis07}. Moreover, they often offer much stronger stability and robustness properties than an EKF and are simpler to tune and implement.

In fact, most existing attitude observers are based on the assumption of weak accelerations of the vehicle so that the gravity direction estimate can be approximated by the accelerometer measurements \cite{tmrm07,mhp08,ms10cep,hua14}. However, the accuracy of the estimated attitude using this assumption is far from satisfactory when the vehicle undergoes sustained accelerations. To get rid of this unsatisfactory assumption some authors have considered attitude observers ``aided'' by complementary linear velocity measurements in the inertial frame \cite{ms08IFAC,hua10cep,rt11CDC} or body-fixed frame~\cite{BonnaMR2006ACC,BonnabelITAC08,troniICRA13,Dukan13}. This paper focuses on the latter category. In \cite{Dukan13}, an observer is proposed and tested experimentally on an underwater vehicle, but without convergence and stability analysis; in the same context \cite{troniICRA13} uses the numerical derivative of the linear velocity to recover the gravity direction estimate, which is sensitive to measurement noise. On the other hand, \cite{BonnabelITAC08} proposes a ``general'' invariant observer with several nice geometric properties among which 1) the local exponential stability  around any trajectory of the system, 2) the local decoupling of the roll and pitch estimation from magnetometer measurements.

In this paper, we propose two invariant observers which can be seen as particular cases of the general invariant observer of~\cite{BonnabelITAC08}. The interest is that we can now guarantee the (almost) global asymptotical stability and the global decoupling of the roll and pitch estimation from magnetometer measurements, while ensuring good local convergence properties and an easy tuning.

\section{Preliminary Material} \label{sec:preliminary}

\subsection{Notation}
\noindent$\bullet\,$ $\{e_1,e_2,e_3\}$ denotes the canonical basis of $\mathbb{R}^3$. The notation $(\cdot)_\times$ denotes the skew-symmetric matrix associated with the cross product, i.e., $u_\times v = u \times v, \forall u, v \in \mathbb{R}^3$.\\
$\bullet\,$ Let $\{\mathcal{I}\}$ denote an inertial frame attached to the earth,  typically chosen as the north-east-down (NED) frame. Let $\{\mathcal B\}$ be a body-fixed frame attached to the vehicle.\\
$\bullet\,$ Let $v \in \mathbb{R}^3$ denote the vehicle's linear velocity, expressed in $\{\mathcal B\}$. Let $\omega \in \mathbb{R}^3$ denote the angular velocity, expressed in $\{\mathcal B\}$, of the frame $\{\mathcal B\}$ with respect to the frame $\{\mathcal I\}$.\\
$\bullet\,$ The vehicle's attitude is represented by a rotation matrix $R \in \text{SO(3)}$ of the frame $\{\mathcal B\}$ relative to  $\{\mathcal I\}$. Let $\phi$, $\theta$ and $\psi$ denote the roll, pitch and yaw Euler angles. By representing the gravity direction by the unit vector $\gamma \triangleq R^T e_3 = [ -\sin\theta, \ \sin\phi \cos\theta, \ \cos\phi \cos\theta]^T$, one deduces that roll and pitch Euler angles can be (locally) uniquely determined from $\gamma$, except singularities corresponding to $\theta = \pm \pi/2$.

\subsection{System Equations and Measurements}
The attitude satisfies the differential equation \vspace{-0.1cm}
\begin{equation}\label{attitudeSystem}
\dot R = R \omega_\times, \vspace{-0.1cm}
\end{equation}
from which the dynamics of $\gamma$ ($=R^Te_3$) are deduced as\vspace{-0.1cm}
\begin{equation}\label{attitudegammaSystem}
\dot \gamma = \gamma \times \omega. \vspace{-0.1cm}
\end{equation}
Assume that the vehicle is equipped with a velocity sensor to measure $v$.
Additionally, it is also equipped with an Inertial Measurement Unit (IMU) consisting of a 3-axis gyrometer and a 3-axis accelerometer. The gyrometer and
accelerometer respectively provide the measurement of the angular velocity $\omega$ and {\it specific acceleration} $a_{\mathcal B} \in \mathbb{R}^3$, expressed in $\{\mathcal B\}$. Using the flat non-rotating Earth assumption, one has \cite{BonnabelITAC08}\vspace{-0.2cm}
\begin{equation}\label{accMeasurement}
a_{\mathcal B} = R^T(\ddot x - ge_3) = \dot v - v \times \omega - g \gamma, \vspace{-0.2cm}
\end{equation}
where the vehicle's acceleration and the gravitational acceleration, both expressed in the inertial frame, are $\ddot x \in \mathbb{R}^3$ and $ge_3$. From \eqref{accMeasurement}, it is also convenient to write \vspace{-0.1cm}
\begin{equation}\label{dotv}
\begin{array}{lcl}
\dot v &= & v \times \omega + a_{\mathcal B} + g R^T e_3=  v \times \omega + a_{\mathcal B} + g\gamma.
\end{array} \vspace{-0.1cm}
\end{equation}
In many IMUs, a 3-axis magnetometer is also integrated to provide the measurement of the Earth's magnetic field vector $m_{\mathcal B} \in \mathbb{R}^3$, expressed in the body-fixed frame. One verifies that $m_{\mathcal B} = R^T m_{\mathcal I}$, with $m_{\mathcal I}$ the Earth's magnetic field vector expressed in the inertial frame. It is reasonable to assume that the Earth's magnetic field vector and the gravity direction are non-collinear, i.e. $m_{\mathcal I} \times e_3 \neq 0$.

In summary, the observer design presented in the next section will be based on the following system \vspace{-0.1cm}
\begin{equation}\label{sys:systemComplet}
\left\{
\begin{array}{lcl}
\dot v &=& v \times \omega + a_{\mathcal{B}} + g R^T e_3  \\
\dot R &=& R \omega_\times
\end{array}\right. \vspace{-0.1cm}
\end{equation}
using $(v, \omega, a_{\cal B}, m_{\cal B})$ as measurements.

\section{Observer Design} \label{sec:observerDesign}

\subsection{Invariant Observer Design}
Let $\hat v \in \mathbb{R}^3$ and $\hat R \in \text{SO(3)}$ denote the estimates of $v$ and $R$, respectively. Consider the nonlinear observer system \vspace{-0.1cm}
\begin{equation}\label{sys:observerGeneral}
\left\{
\begin{array}{lcl}
\dot {\hat v} &=& \hat v \times\omega + a_{\mathcal{B}} + g \hat R^T e_3 + \sigma_v\\
\dot {\hat R} &=& \hat R (\omega + \sigma_{\!R})_\times
\end{array} \right.\vspace{-0.1cm}
\end{equation}
where $\sigma_v \in \mathbb{R}^3$ and $\sigma_{\!R} \in \mathbb{R}^3$ are the innovation terms to be designed in order to ensure the following objectives:
\begin{enumerate}[O.1)]
\item the convergence of $(\hat v, \hat R)$ to $(v,R)$ and the stability of the equilibrium $(\hat v, \hat R)=(v,R)$;
\item the decoupling of roll and pitch estimation from magnetometer measurements.
\end{enumerate}

A ``general'' invariant observer in the form of quaternions has been proposed in \cite{BonnabelITAC08}, which is equivalently rewritten as \vspace{-0.2cm}
\begin{equation}\label{sys:observerInvariantGeneral}
\left\{
\begin{array}{lcl}
\!\!\dot {\hat v} &\!\!\!=\!\!\!& \hat v \times\omega + a_{\mathcal{B}} + g \hat R^T e_3 + \hat R^T ({\cal L}_v^v e_v + {\cal L}_m^v e_m)\\
\!\!\dot {\hat R} &\!\!\!=\!\!\!& \hat R \omega_\times + ({\cal L}_v^R e_v + {\cal L}_m^R e_m)_\times \hat R
\end{array} \right. \!\!\! \vspace{-0.1cm}
\end{equation}
where ${\cal L}_v^v$, ${\cal L}_m^v$, ${\cal L}_v^R$, and ${\cal L}_m^R$ are $3\times 3$ gain matrices whose entries may depend on the invariant errors
$e_v \triangleq \hat R(v-\hat v)$ and $e_m \triangleq m_{\cal I} - \hat R m_{\cal B}$.
Motivated by the fact that the convergence and stability and the decoupling results proved in \cite{BonnabelITAC08} are only {\em local}, we propose two following observers which can be seen as particular cases of \eqref{sys:observerInvariantGeneral}, but for which {\em almost global} convergence and stability and {\em global} decoupling properties can be established hereafter:

\vspace{0.1cm}
\noindent $\bullet\,$ {\em Observer 1:} This observer is given by \eqref{sys:observerGeneral}, with $\sigma_v$ and $\sigma_R$ defined by \vspace{-0.1cm}
\begin{equation}\label{innovationObs1}
\left\{\begin{array}{lcl}
\sigma_v&\!\!\!\!\triangleq\!\!\!\!& k_1^v \tilde v  - k_2^v \hat \gamma \times (\hat \gamma \times \tilde v)\\
\sigma_{\! R}&\!\!\!\!\triangleq\!\!\!\!& k_1^r \tilde v \times \hat \gamma
+ k_2^{r} ((m_{\cal B} \!\times\! \hat m_{\cal B})^T\hat \gamma) \hat\gamma
\end{array}\right.\!\!\! \vspace{-0.cm}
\end{equation}
where $k_1^v,\, k_2^v,\, k_1^r,\, k_2^r$ are positive constant gains
and \vspace{-0.cm}
\begin{equation}\label{barvhatgamhatmB}
\tilde v \triangleq v - \hat v, \quad \hat \gamma \triangleq \hat R^T e_3, \quad \hat m_{\cal B} \triangleq \hat R^Tm_{\cal I} \vspace{-0.cm}
\end{equation}

\noindent $\bullet\,$ {\em Observer 2:} This observer is also given by \eqref{sys:observerGeneral}, but with $\sigma_v$ and $\sigma_R$ defined by \vspace{-0.1cm}
\begin{equation}\label{innovationObs2}
\left\{\begin{array}{lcl}
\sigma_v&\!\!\!\!\triangleq\!\!\!\!& k_1^v \tilde v  - k_2^v \hat \gamma \times \!(\hat \gamma\! \times \!\tilde v)\! - k_{1}^r \tilde v \times\! (\tilde v \times\! \hat \gamma)\\
\sigma_{\!R}&\!\!\!\!\triangleq\!\!\!\!& k_1^r \tilde v \times \hat \gamma
+ k_2^{r} ((m_{\cal B} \times \hat m_{\cal B})^T\hat \gamma) \hat\gamma
\end{array}\right.\!\!\! \vspace{-0.cm}
\end{equation}
where $k_1^v,\, k_2^v,\, k_1^r,\, k_2^r$ are positive constant gains, and $\tilde v$, $\hat \gamma$, $\hat m_{\cal B}$ are defined by \eqref{barvhatgamhatmB}.

The sole difference between the two observers is the ``quadratic'' term $-k_{1}^r \tilde v \!\times\! (\tilde v \!\times \!\hat \gamma)$ involved in the definition \eqref{innovationObs2} of $\sigma_v$ of Observer 2. We will explain later why this term has been introduced.

\vspace{-0.3cm}
\begin{lem} \label{lemma:ErrorDyn}
The dynamics of the invariant state error $(\bar v, \bar R)$ defined by \vspace{-0.1cm}
\begin{equation} \label{error}
 \bar v \triangleq R (v-\hat v)  \quad \bar R \triangleq R \hat R^T \vspace{-0.1cm}
\end{equation}
for both Observers 1 and 2 are autonomous.
\end{lem}
\vspace{-0.1cm}
\begin{proof} The proof is similar to the one in \cite{BonnabelITAC08}. Using \eqref{sys:systemComplet}, \eqref{sys:observerGeneral} and \eqref{error}, one easily deduces

\vspace{-0.3cm}
\begin{equation} \label{sys:errorGen}
\left\{
\begin{array}{lcl}
\dot{\bar v} &=&  g(I - \bar R)e_3 - \bar \sigma_v \\
\dot{\bar R} &=&  -\bar \sigma_{\!R\times}  \bar R
\end{array} \right. \vspace{-0.1cm}
\end{equation}
with $\bar \sigma_v \triangleq R \sigma_v, \, \bar \sigma_{\!R} \triangleq R \sigma_{\!R}$. When $\sigma_v$ and $\sigma_{\! R}$ are given by \eqref{innovationObs1} for Observer 1, using the identity $R (u \times v) = (R u) \times (Rv)$, $\forall u,v \in \mathbb{R}^3$, $\forall R \in \text{SO(3)}$, one deduces \vspace{-0.1cm}
\begin{IEEEeqnarray}{RCL}
  \bar \sigma_v &=& k_1^v \bar v - k_2^v (\bar R e_3) \times ((\bar R e_3) \times\bar v) \IEEEyessubnumber \label{eq:RsigmavObs1} \\
  \bar \sigma_{\! R} &=& k_1^r \bar v \times \bar R e_3+ k_2^r ((m_{\cal I} \times \bar R m_{\cal I})^T \bar R e_3) \bar R e_3
  \IEEEyessubnumber \label{eq:RsigmaRObs1} \vspace{-0.1cm}
\end{IEEEeqnarray}
When $\sigma_v$ and $\sigma_{\! R}$ are given by \eqref{innovationObs2} for Observer 2, one verifies that $\bar \sigma_{\! R}$ still satisfies \eqref{eq:RsigmaRObs1}, while $\bar\sigma_v$ is given by \vspace{-0.2cm}
\begin{equation} \label{eq:RsigmavObs2}
\begin{array}{lcl}
\!\!\bar \sigma_v &\!\!\!\!\!=\!\!\!\!\!& k_1^v \bar v \!-\! k_2^v (\bar R e_3) \!\!\times\! ((\bar R e_3)\! \!\times\!\bar v) \!-\! k_1^r  \bar v\!\times\! (\bar v \!\times\! \bar R e_3)
\end{array} \!\! \vspace{-0.1cm}
\end{equation}
From here, the conclusion is straightforward.
\end{proof}

\vspace{-0.1cm}
\subsection{Reduction to Gravity Direction Estimation}
The gravity direction expressed in the body-fixed frame can be represented by the vector $\gamma\, (=\!R^T e_3 \!\in\! S^2)$. Its estimate $\hat \gamma \in S^2$ can be calculated from the estimate $\hat R\in \text{SO(3)}$ provided by Observers 1 or 2 as $\hat \gamma = \hat R^T e_3$. It can also be obtained from an observer directly designed on $\mathbb{R}^3 \times S^2$ as a result of the following lemma.

\vspace{-0.3cm}
\begin{lem} \label{lemma:ObserverGamma}
Observers 1 and 2 can be reduced to the following observers of $v$ and $\gamma$ (that we term ``$\gamma$--observers''):
\begin{itemize}
\item {$\gamma$-Observer 1:} \vspace{-0.1cm}
\begin{equation}\label{sys:GammaObs1}
\left\{
\begin{array}{lcl}
\dot {\hat v} &=& \hat v \times\omega + a_{\mathcal{B}} + g \hat \gamma + \sigma_v^\gamma\\
\dot {\hat \gamma} &=& \hat \gamma \times (\omega + \sigma_{\!R}^\gamma) \\
\sigma_v^\gamma &\triangleq& k_1^v (v -\hat v)  - k_2^v \hat \gamma \times (\hat \gamma \times (v-\hat v)) \\
\sigma_{\! R}^\gamma &\triangleq& k_1^r (v-\hat v) \times \hat \gamma
\end{array} \right.
\end{equation}
where $k_1^v$, $k_2^v$, $k_1^r$ are positive constant gains.
\item {$\gamma$-Observer 2:} \vspace{-0.1cm}
\begin{equation}\label{sys:GammaObs2}
\left\{
\begin{array}{lcl}
\dot {\hat v} &=& \hat v \times\omega + a_{\mathcal{B}} + g \hat \gamma + \sigma_v^\gamma\\
\dot {\hat \gamma} &=& \hat \gamma \times (\omega + \sigma_{\!R}^\gamma) \\
\sigma_v^\gamma &\triangleq& k_1^v (v -\hat v)  - k_2^v \hat \gamma \times (\hat \gamma \times (v-\hat v)) \\
&& - k_1^r (v-\hat v) \times ((v-\hat v) \times \hat\gamma)\\
\sigma_{\! R}^\gamma &\triangleq& k_1^r (v-\hat v) \times \hat \gamma
\end{array} \right.
\end{equation}
where $k_1^v$, $k_2^v$, $k_1^r$ are positive constant gains.
\end{itemize}
In addition, these two $\gamma$--observers are independent of magnetometer measurements (i.e., $m_{\cal B}$).
\end{lem}
\vspace{-0.2cm}
\begin{proof}
The expression of $\dot{\hat v}$ in \eqref{sys:GammaObs1} (resp. \eqref{sys:GammaObs2}) is straightforwardly obtained from \eqref{sys:observerGeneral} and \eqref{innovationObs1} (resp. \eqref{innovationObs2}) by replacing $\hat R^T e_3$ by $\hat \gamma$. As for the dynamics of $\hat \gamma$, by differentiating $\hat\gamma=\hat R^T e_3$ and using the expression of $\dot {\hat R}$ given by \eqref{sys:observerGeneral}--\eqref{innovationObs1} (or \eqref{sys:observerGeneral}--\eqref{innovationObs2}) one deduces \vspace{-0.1cm}
\[
\dot{\hat \gamma} = \hat \gamma \times (\omega + \sigma_{\!R}) = \hat \gamma \times (\omega + \sigma_{\!R}^\gamma) \vspace{-0.1cm}
\]
where the latter equality is obtained using \vspace{-0.1cm}
\[
\hat\gamma \times [k_2^{r} ((m_{\cal B} \!\times\! \hat m_{\cal B})^T\hat \gamma) \hat\gamma
] = 0 \vspace{-0.1cm}
\]
Finally, the statement about the independence of the $\gamma$--observers \eqref{sys:GammaObs1} and \eqref{sys:GammaObs2} on $m_{\cal B}$ is straightforward.
\end{proof}

\vspace{0.1cm}
The latter statement in Lemma \ref{lemma:ObserverGamma} implies that the objective O.2 is guaranteed globally. With respect to the statement in Lemma \ref{lemma:ErrorDyn} that the dynamics of the invariant estimation state errors $(\bar v,\bar R)$ are autonomous, a similar result is now given.

\vspace{-0.3cm}
\begin{lem} \label{lemma:ErrorDynGamma}
The dynamics of the invariant state errors $(\bar v, \bar \gamma)$ defined by \vspace{-0.2cm}
\begin{equation} \label{errorGamma}
\bar v \triangleq R (v-\hat v), \quad \bar \gamma \triangleq R \hat \gamma = \bar R e_3 \vspace{-0.2cm}
\end{equation}
for both $\gamma$--Observers 1 and 2, given by \eqref{sys:GammaObs1} and \eqref{sys:GammaObs2}, respectively, are autonomous.
\end{lem}
\vspace{-0.2cm}
\begin{proof}
For both $\gamma$--Observers 1 and 2, one verifies from \eqref{attitudegammaSystem}, \eqref{dotv}, \eqref{sys:GammaObs1} (or \eqref{sys:GammaObs2}) that \vspace{-0.2cm}
\[
\begin{array}{lcl}
\dot{\bar v} &\!\!\!\!= \!\!\!\!& R \omega_\times (v -\hat v) + R (\dot v-\dot {\hat v})= g (e_3 - \bar\gamma) - R \sigma_v^\gamma \vspace{-0.1cm}
\end{array}
\]
\[
\dot{\bar \gamma} = R \omega_\times \hat\gamma + R ( \hat\gamma \times (\omega + \sigma_{\! R}^\gamma)) =
\bar \gamma \times (R \sigma_{\! R}^\gamma) = -k_1^r \bar \gamma \times (\bar \gamma \times \bar v) \vspace{-0.1cm}
\]
Consequently, for $\gamma$--Observer 1 (i.e., \eqref{sys:GammaObs1}) one obtains \vspace{-0.1cm}
\begin{equation} \label{sys:errorGammaObs1}
\left\{
\begin{array}{lcl}
\dot{\bar v} &\!\!\!\!= \!\!\!\! &  g (e_3 - \bar \gamma) - k_1^v  \bar v + k_2^v \bar\gamma\times (\bar \gamma \times \bar v) \\
\dot{\bar \gamma} &\!\!\!\!= \!\!\!\!& -k_1^r \bar \gamma \times (\bar \gamma \times \bar v)
\end{array} \right.
\end{equation}
On the other hand, for $\gamma$--Observer 2 (i.e., \eqref{sys:GammaObs2}) it yields
\begin{equation} \label{sys:errorGammaObs2}
\left\{
\begin{array}{lcl}
\!\!\!\dot{\bar v} &\!\!\!\!= \!\!\!\!&  g (e_3 \!-\! \bar \gamma)
- k_1^v  \bar v \!+\! k_2^v \bar\gamma\!\times\! (\bar \gamma \!\times\! \bar v)\!+\! k_1^r \bar v \!\times\! (\bar v \!\times\!\bar\gamma) \\
\!\!\!\dot{\bar \gamma} &\!\!\!\!= \!\!\!\!& -k_1^r \bar \gamma \times (\bar \gamma \times \bar v)
\end{array} \right.\!\!\!\!
\end{equation}
The conclusion then directly follows.
\end{proof}

In the following, convergence and stability analyses of the error systems \eqref{sys:errorGammaObs1} and \eqref{sys:errorGammaObs2} are provided.

\vspace{-0.3cm}
\begin{propo} \label{propo:gammaObs1} ($\gamma$--Observer 1)--
Consider the autonomous error dynamics \eqref{sys:errorGammaObs1} and assume that the observer gains $k_1^v$, $k_2^v$, $k_1^r$ are chosen positive and satisfying the following condition: \vspace{-0.1cm}
\begin{equation}\label{eq:gainCond}
k_1^r \leq \frac{k_1^v k_2^v}{g}
\vspace{-0.1cm}
\end{equation}
Then, the following properties hold:
\begin{enumerate}
\item System \eqref{sys:errorGammaObs1} has only two isolated equilibrium points $(\bar v,\bar \gamma) = (0,e_3)$ and $(\bar v,\bar \gamma) = (\frac{2g}{k_1^v} e_3,-e_3)$. For all initial condition $(\bar v(0),\bar \gamma(0))$, the error trajectory $(\bar v(t),\bar \gamma(t))$ converges to one of these two equilibria.
\item The equilibrium $(\bar v,\bar \gamma) = (0,e_3)$ is almost-globally asymptotically stable and locally exponentially stable.
\item The equilibrium $(\bar v,\bar \gamma) = (\frac{2g}{k_1^v} e_3,-e_3)$ is unstable.
\end{enumerate}
\end{propo}
\vspace{-0.1cm}
\begin{proof}
Consider the Lyapunov function candidate \vspace{-0.1cm}
\begin{equation}\label{eq:lyapunovL0}
\begin{split}
{\cal L}_0 &\triangleq \frac{1}{2}|\bar v|^2 + \frac{gk_2^v}{2 k_1^v k_1^r} |e_3 - \bar \gamma|^2 - \frac{g}{k_1^v} \bar v^T (e_3 -\bar\gamma)
\end{split}
\end{equation}

\vspace{-0.cm}
\noindent which is positive-definite under condition \eqref{eq:gainCond}.
From \eqref{sys:errorGammaObs1} and \eqref{eq:lyapunovL0}, one verifies that the time-derivative of ${\cal L}_0$ satisfies \vspace{-0.2cm}
\begin{equation*}\label{eq:lyapunovdotL0}
\begin{array}{lcl}
\!\!\!\dot{\cal L}_0
%&\!\!\!\!\!= \!\!\!\!\!&- k_1^v  |\bar v|^2 - \frac{g^2}{k_1^v} |e_3-\bar\gamma|^2 + 2 g \bar v^T(e_3 - \bar \gamma) \\
%&& + (k_2^v - \frac{g k_1^r}{k_1^v}) \bar v^T(\bar\gamma\times (\bar \gamma \times \bar v)) \\
&\!\!\!\!\!= \!\!\!\!\!&-k_1^v |\bar v - \frac{g}{k_1^v} (e_3 - \bar\gamma)|^2 -(k_2^v -\frac{gk_1^r}{k_1^v}) |\bar \gamma \times \bar v|^2
\end{array}\!\! \vspace{-0.1cm}
\end{equation*}
which is negative-(semi)definite under condition \eqref{eq:gainCond}. Since System \eqref{sys:errorGammaObs1} is autonomous, by application of LaSalle's theorem, one deduces the convergence of $\dot{\cal L}_0$ to zero. This in turn implies the convergence of $\delta \triangleq \bar v - \frac{g}{k_1^v} (e_3 - \bar\gamma)$ to zero, and additionally the convergence of $\bar \gamma \times \bar v$ to zero if $k_1^r  < \frac{k_1^v k_2^v}{g}$. Contrarily, if $k_1^r  = \frac{k_1^v k_2^v}{g}$, only the convergence of $\delta$ to zero can be deduced. Let us now prove the convergence of $e_3\times\bar\gamma$ to zero for the two possible cases satisfying \eqref{eq:gainCond}.

$\bullet$ Case $k_1^r  < \frac{k_1^v k_2^v}{g}$: The convergence of $\bar v$ to $\frac{g}{k_1^v} (e_3 - \bar\gamma)$ can be deduced from the definition of $\delta$ and its convergence to zero (proved previously). This implies that
$\bar v \times \bar \gamma$ converges to $\frac{g}{k_1^v}  e_3 \times \bar \gamma$, which must converge to zero since $\bar v \times \bar \gamma$ converges to zero (proved previously).

$\bullet$ Case $k_1^r  = \frac{k_1^v k_2^v}{g}$: It is easily deduced from \eqref{sys:errorGammaObs1} and the definition of $\delta$ that $\dot \delta = -k_0 \delta$. On the other hand, the definition of $\delta$ implies $\bar v =\delta + \frac{g}{k_1^v} (e_3 - \bar\gamma)$. Then, it can be verified  from \eqref{sys:errorGammaObs1} that \vspace{-0.1cm}
\begin{equation} \label{dotTerm}
\begin{array}{lcl}
\frac{d}{dt}(1\!-\! e_3^T \bar\gamma ) &\!\!\!\!\!\! \!\! \!\! =\!\!\!\!\!\! \!\! \!\! & -k_1^r e_3^T (\bar\gamma \times (\bar v\times\bar\gamma))  \\
&\!\! \!\!\!\!\!\!\!  =\!\! \!\!\!\!\!\!\!  & - k_1^r ((\delta + \frac{g}{k_1^v} (e_3 - \bar\gamma))\times \bar\gamma)^T (e_3 \times \bar\gamma)\\
&\!\! \!\!\!\!\!\!\!  =\!\! \!\!\!\!\!\!\!  & -k_1^r (\delta\times\bar\gamma)^T (e_3 \times \bar\gamma) -\frac{gk_1^r}{k_1^v} |e_3 \times \bar\gamma|^2 \\
&\leq& k_1^r |\delta| |e_3 \times \bar\gamma| -\frac{gk_1^r}{k_1^v} |e_3 \times \bar\gamma|^2
\end{array} \vspace{-0.2cm}
\end{equation}
Now, consider the following Lyapunov function candidate \vspace{-0.1cm}
\begin{equation}\label{eq:lyapunovL1}
{\cal L}_1 \triangleq \frac{1}{2} |\delta|^2 + \frac{2g}{k_1^r}(1- e_3^T \bar\gamma)  \vspace{-0.1cm}
\end{equation}
One deduces from \eqref{dotTerm}, \eqref{eq:lyapunovL1} and $\dot \delta = -k_0 \delta$ that \vspace{-0.cm}
\begin{equation*}\label{eq:lyapunovdotL1}
\dot {\cal L}_1 \leq - k_1^v |\delta|^2 - ({2g^2}/{k_1^v}) |e_3 \times \bar\gamma|^2 +  2g |\delta| |e_3 \times \bar\gamma| \vspace{-0.cm}
\end{equation*}
Since $k_1^v |\delta|^2 + ({2g^2}/{k_1^v}) |e_3 \times \bar\gamma|^2 \geq 2 \sqrt{2} g |\delta| |e_3 \times \bar\gamma|$ using Young's inequality, there exist two positive numbers $\alpha_1$ and $\alpha_2$ such that
$\dot {\cal L}_1 \leq - \alpha_1 |\delta|^2 - \alpha_2 |e_3 \times \bar\gamma|^2$.
Then, by application of LaSalle's theorem, one deduces the convergence of $\dot {\cal L}_1$ to zero, which implies that $\delta$ and $e_3 \times \bar\gamma$ also converge to zero.

Therefore, for both cases $\bar \gamma$ converges to either $e_3$ (desired) or $-e_3$ (undesired), which in turn implies that $\bar v$ converges to either zero (desired) or $\frac{2g}{k_1^v}e_3$ (undesired). The stability is directly deduced from the expressions of either $({\cal L}_0, \dot{\cal L}_0)$ or $({\cal L}_1, \dot{\cal L}_1)$.

We continue to prove that the ``desired'' equilibrium $(\bar v,\bar \gamma) =(0,e_3)$ is locally exponentially stable and the ``undesired'' one $(\bar v,\bar \gamma) =(({2g}/{k_1^v})e_3,-e_3)$ is unstable. To this purpose, it suffices to study the stability of the linearized system about each equilibrium.
For instance, the linearized system about the equilibrium $(\bar v,\bar \gamma) =(0,e_3)$ satisfies \vspace{-0.cm}
\[\left\{
\begin{array}{lcl}
\dot{\bar v} &=& g(e_3-\bar\gamma) -k_1^v \bar v
+ k_2^v e_3 \times(e_3 \times \bar v) \\
\dot{\bar\gamma} &=&  - k_1^r e_3 \times (e_3\times \bar v)
\end{array}\right. \vspace{-0.cm}
\]
which can be decomposed into three subsystems \vspace{-0.1cm}
\begin{IEEEeqnarray}{RCL}
\begin{bmatrix} \dot {\bar v}_{i} \\ \dot{\bar\gamma}_{i}\end{bmatrix} &=& \begin{bmatrix}-(k_1^v +k_2^v) & -g \\ k_1^r & 0 \end{bmatrix}
\begin{bmatrix} {\bar v}_{i} \\  {\bar\gamma}_{i}\end{bmatrix}, \,i=1,2 \IEEEyessubnumber\label{eq:sysLinGoodEqui-1} \\
\dot{\bar v}_3 &=& -k_1^v \bar v_3  \IEEEyessubnumber\label{eq:sysLinGoodEqui-3} \vspace{-0.cm}
\end{IEEEeqnarray}
By application of Hurwitz criteria, one easily deduces that the origin of these three subsystems is stable for any set of positive constant gains $(k_1^v,k_2^v,k_1^r)$. Similarly, one can easily prove that the remaining ``undesired'' equilibrium $(\bar v,\bar \gamma) =(({2g}/{k_1^v})e_3,-e_3)$ is unstable by analysing the linearized system about this equilibrium. Therefore, the equilibrium $(\bar v,\bar \gamma) =(0,e_3)$ is almost globally asymptotically stable and locally exponentially stable.
\end{proof}

\vspace{0.2cm}
\noindent {\bf Remarks:}

1) The gain condition \eqref{eq:gainCond} is only sufficient for the almost global asymptotical stability of the equilibrium  $(\bar v,\bar \gamma) =(0,e_3)$. If \eqref{eq:gainCond} is not satisfied, in view of the linearized system \eqref{eq:sysLinGoodEqui-1}--\eqref{eq:sysLinGoodEqui-3} one still deduces the local exponential stability of the equilibrium  $(\bar v,\bar \gamma) =(0,e_3)$, for any set of positive constant gains $(k_1^v,k_2^v,k_1^r)$.

2) The linearized subsystems \eqref{eq:sysLinGoodEqui-1}, with $i=1,2$, have identical form, with characteristic polynomial given by
$P(\lambda) =\lambda^2 + (k_1^v + k_2^v)\lambda + gk_1^r$. Using Young's inequality and condition \eqref{eq:gainCond}, one verifies that the determinant of $P(\lambda)$ is positive, i.e., $\Delta_P = (k_1^v + k_2^v)^2 - 4gk_1^r \geq 4 k_1^vk_2^v - 4gk_1^r > 0$. This implies that $P(\lambda)$ can only possess two negative real poles. This in turn implies that we cannot impose imaginary poles for these subsystems while respecting condition \eqref{eq:gainCond}. This limitation is usually not critical in practice. However, from a theoretical standpoint we still want to obtain a ``stronger'' result in the sense that the almost global asymptotical stability property of $\gamma$--Observer 1 is still ensured, while the poles of the linearized system about the ``desired'' equilibrium can be arbitrarily chosen (with negative real part). Such a motivation has led us to introduce the ``quadratic'' term $-k_1^r \tilde v \times (\tilde v \times \hat\gamma)$ in the innovation term $\sigma_v$ of Observer 1, yielding Observer 2. This then yields the error system \eqref{sys:errorGammaObs2} of $\gamma$--Observer 2, which only differs from \eqref{sys:errorGammaObs1} by the ``quadratic'' term $-k_1^r \bar v \times (\bar v \times \bar\gamma)$. Since this term is neglected in the linearized system of \eqref{sys:errorGammaObs2} about the equilibrium $(\bar v,\bar \gamma) =(0,e_3)$, the linearized systems of both \eqref{sys:errorGammaObs1} and \eqref{sys:errorGammaObs2} are identical. As shown in Proposition~\ref{propo:gammaObs2} (presented below), condition \eqref{eq:gainCond} is no longer required for the almost-global asymptotical stability of the error system  \eqref{sys:errorGammaObs2}, which in turn implies that more freedom for the choice of poles is available for gain tuning. This is an advantage of Observer 2 (resp. $\gamma$--Observer 2) with respect to Observer 1 (resp. $\gamma$--Observer 1). Conversely, introducing an addition term in the observer may make it more sensitive to measurement noises. This means that both observers have advantages with respect to each other.

\vspace{-0.3cm}
\begin{propo} \label{propo:gammaObs2} ($\gamma$--Observer 2)--
Consider the autonomous error dynamics \eqref{sys:errorGammaObs2}, with $k_1^v$, $k_2^v$, $k_1^r$ chosen positive.
Then, all  properties given in Proposition \ref{propo:gammaObs1} hold.
\end{propo}
\vspace{-0.2cm}
\begin{proof} Using \eqref{sys:errorGammaObs2} and the identity $u \times (v \times w) = v (u^T w) - w (u^Tv)$, $\forall u,v,w \in \mathbb{R}^3$, one verifies that \vspace{-0.1cm}
\begin{IEEEeqnarray}{RCL}
\!\!\frac{d}{dt} (\bar v\! \times\! \bar \gamma) &=& g \,e_3 \times \bar\gamma - (k_1^v + k_2^v)\, \bar v \times \bar\gamma
\label{dotvbartimesgammabar} \vspace{-0.1cm}
\end{IEEEeqnarray}
From \eqref{sys:errorGammaObs2} and \eqref{dotvbartimesgammabar}, it is clear that the dynamics of $(\bar v \times \bar \gamma, \bar\gamma)$ are autonomous. Consider the following positive function \vspace{-0.1cm}
\begin{equation}\label{eq:lyapunovS0}
{\cal S}_0 \triangleq {1}/{2} |\bar v \times \bar \gamma|^2 + {g}/{k_1^r} (1 - e_3^T \bar \gamma) \vspace{-0.cm}
\end{equation}
Using \eqref{sys:errorGammaObs2}, \eqref{dotvbartimesgammabar} and \eqref{eq:lyapunovS0}, one deduces \vspace{-0.1cm}
\begin{IEEEeqnarray}{RCL}
\dot{\cal S}_0  &=& - (k_1^v + k_2^v)\, |\bar v \times \bar\gamma|^2 \label{eq:dotlyapunovS0} \vspace{-0.1cm}
\end{IEEEeqnarray}
From here, by application of LaSalle's theorem, one deduces the convergence of $\dot{\cal S}_0$ and, subsequently, of $\bar v \times \bar\gamma$ to zero.

From \eqref{eq:lyapunovS0} and \eqref{eq:dotlyapunovS0}, one deduces the boundedness of $\bar v \times \bar\gamma$. One then verifies that $\dot{\bar \gamma}$ is also bounded, which implies the uniform continuity of $\bar\gamma$. Then, by application of the extended Barbalat's lemma (see, e.g., \cite{micaelli1993}) to \eqref{dotvbartimesgammabar} one deduces the convergence of $\frac{d}{dt} (\bar v\! \times\! \bar \gamma)$ to zero. This in turn implies the convergence of $e_3 \times \bar\gamma$ to zero. Therefore, $\bar \gamma$ converges to either $e_3$ or $-e_3$.

Since $\bar v \!\times\! \bar\gamma$ converges to zero, the zero-dynamics of $\dot{\bar v}$ are \vspace{-0.1cm}
\begin{equation}\label{dotvbarZerodyn}
\dot{\bar v} =  g (e_3 - \bar \gamma) - k_1^v  \bar v \vspace{-0.1cm}
\end{equation}
Then, the convergence of $\bar \gamma$ to either $e_3$ or $-e_3$ associated with the zero-dynamics \eqref{dotvbarZerodyn} ensures the convergence of $\bar v$ to either zero or $(2g/k_1^v)e_3$.

Finally, the proof of almost-global asymptotical stability and local exponential stability of the ``desired'' equilibrium $(\bar v,\bar \gamma)=(0,e_3)$ and the proof of instability of the ``undesired'' one $(\bar v,\bar \gamma) =(({2g}/{k_1^v})e_3,-e_3)$ proceed analogously to the proof of Proposition \ref{propo:gammaObs1}.
\end{proof}

\subsection{Stability Analysis for Observers 1 and 2}

In order to analyze the asymptotic stability of the observer trajectory of Observers 1 and 2 to the system trajectory, it is more convenient to consider the dynamics of the state errors $(\bar v, \bar R)$ defined by \eqref{error} and prove that their trajectory converges to $(0,I)$, with $I$ the identity element of $\text{SO(3)}$.

\vspace{-0.3cm}
\begin{theo} \label{theo:Observer1} (Observer 1)-- Consider System \eqref{sys:systemComplet} and Observer 1 (i.e. \eqref{sys:observerGeneral}+\eqref{innovationObs1}). Assume that condition \eqref{eq:gainCond} for the observer gains $k_1^v, k_2^v, k_1^r$ is satisfied and $m_{\cal I} \times e_3 \neq 0$. Then, the following properties hold:
\begin{enumerate}
\item The dynamics of the estimate errors $(\bar v,\bar R)$ have only four isolated equilibria, one of which is $(\bar v,\bar R) = (0,I)$.
\item The equilibrium $(\bar v,\bar R) = (0,I)$ is locally exponentially stable and almost globally asymptotically stable; and the other three equilibria of $(\bar v,\bar R)$ are unstable. Thus, for almost all initial conditions $(\hat v (0), \hat R(0))$, the trajectory $(\hat v(t), \hat R(t))$ converges to the system trajectory $(v(t), R(t))$.
\item The dynamics of $(\hat v, \hat R^T e_3)$ are independent of $m_{\cal B}$.
\end{enumerate}
\end{theo}
\vspace{-0.2cm}
\begin{proof} Property 3 is a direct result of Lemma \ref{lemma:ObserverGamma}. We now prove Property 1.
First, recall that the dynamics of $(\bar v,\bar R)$ are given by \eqref{sys:errorGen}, with $\bar\sigma_v$ and $\bar\sigma_{\!R}$ given by \eqref{eq:RsigmavObs1} and \eqref{eq:RsigmaRObs1}, respectively.
As a result of Proposition \ref{propo:gammaObs1}, one ensures the convergence of $(\bar R e_3, \bar v)$ to either $(e_3, 0)$ or $(-e_3, \frac{2g}{k_1^v}e_3)$. For both cases, the term $\bar\sigma_{\!R}$ given by \eqref{eq:RsigmaRObs1} converges exponentially to  \vspace{-0.1cm}
\begin{IEEEeqnarray}{RCL}
\bar\sigma_{\!R} &\rightarrow& k_2^r ((m_{\cal I} \times \bar R m_{\cal I})^T e_3) e_3 \notag\\
&\rightarrow& k_2^r ((\pi_{e_3} m_{\cal I} \times \bar R (\pi_{e_3}m_{\cal I}))^T e_3) e_3 \notag \\
&\rightarrow& k_2^r |\pi_{e_3} m_{\cal I}|^2 ((\bar  m_{\cal I} \times \bar R \bar  m_{\cal I})^T e_3 )e_3 \notag \vspace{-0.1cm}
\end{IEEEeqnarray}
where $\pi_x = I - x x^T$, $\forall x\in \mathbb{R}^3$ denote the projection on the plane orthogonal to $x$, and  $\bar m_{\cal I} \triangleq \frac{\pi_{e_3} m_{\cal I}}{|\pi_{e_3} m_{\cal I}|}$. Consequently, the dynamics of $\bar R$ write \vspace{-0.cm}
\begin{equation}\label{dotRbarZeroDyn}
\dot{\bar R} = - \bar k_2^r ((\bar  m_{\cal I} \times \bar R \bar  m_{\cal I})^T e_3) e_{3\times} \bar R  + \epsilon(\bar v,\bar R)_\times \bar R\vspace{-0.1cm}
\end{equation}
with $\bar k_2^r \triangleq k_2^r |\pi_{e_3} m_{\cal I}|^2$ and a term $\epsilon(\bar v,\bar R)\in \mathbb{R}^3$ remaining bounded and converging exponentially to zero. One can also easily verify that $\bar R$ is uniformly continuous.

Using \eqref{dotRbarZeroDyn}, one verifies that the time-derivative of the positive function ${\cal V} \triangleq 1- \bar  m_{\cal I}^T \bar R \bar  m_{\cal I}$ satisfies $\dot{\cal V} \leq - \bar k_2^r ((\bar  m_{\cal I} \times \bar R \bar  m_{\cal I})^T e_3)^2 + |\epsilon(\bar v,\bar R)|$. Then, by integration one deduces \vspace{-0.1cm}
\[
\int_0^\infty \!\!\!\!\bar k_2^r ((\bar  m_{\cal I} \times \bar R \bar  m_{\cal I})^T e_3)^2 d\tau \leq \int_0^\infty\!\!\!\! |\epsilon(\bar v,\bar R)|d\tau + {\cal V}(0)-{\cal V}(\infty) \vspace{-0.1cm}
\]
From here, one deduces that $\int_0^\infty \!\!((\bar  m_{\cal I} \times \bar R \bar  m_{\cal I})^T e_3)^2 d\tau$ remains bounded since $\cal V$ is bounded and $|\epsilon(\bar v,\bar R)|$ converges exponentially to zero. Then, the application of Barbalat's lemma yields the convergence of $(\bar  m_{\cal I} \times \bar R \bar  m_{\cal I})^T e_3$ to zero.

Now, from the definition of $\bar  m_{\cal I}$, one deduces that this constant unit vector belongs to $\text{Span}(e_1,e_2)$. Thus, there exists a constant angle $\alpha$ such that \vspace{-0.15cm}
\begin{equation*}\label{defRalpha}
\bar  m_{\cal I} = \cos\alpha\, e_1 +\sin\alpha \, e_2
=  {\tiny\begin{bmatrix} \cos\alpha  & -\sin\alpha  & 0 \\ \sin\alpha & \cos\alpha & 0 \\ 0& 0&1 \end{bmatrix}} e_1 = {R}_{\alpha} {e}_1 \vspace{-0.15cm}
\end{equation*}
Since $R_\alpha e_3 =R_\alpha^T e_3 = e_3$, one writes \vspace{-0.15cm}
\[
(\bar  m_{\cal I} \times \bar R \bar  m_{\cal I})^T e_3 = (e_1 \times R_\alpha^T \bar R R_\alpha e_1)^T e_3 = e_2^T(R_\alpha^T \bar R R_\alpha) e_1  \vspace{-0.15cm}
\]
which implies that $e_2^T(R_\alpha^T \bar R R_\alpha) e_1 \rightarrow 0$ using the fact that $(\bar  m_{\cal I} \times \bar R \bar  m_{\cal I})^T e_3 \rightarrow 0$. One also verifies from the convergence $\bar R e_3 \rightarrow \pm e_3$ that $R_\alpha^T \bar R R_\alpha e_3 \rightarrow \pm e_3$. From here, it is straightforward to deduce that $R_\alpha^T \bar R R_\alpha$ converges to one of the four following rotation matrices: $R_1^*\triangleq I$, $R_2^* \triangleq \text{diag}(-1,-1,1)$, $R_3^* \triangleq \text{diag}(-1,1,-1)$, $R_4^* \triangleq \text{diag}(1,-1,-1)$,
where the first two matrices correspond to the case $\bar R e_3 \rightarrow  e_3$, and the last two correspond to the case $\bar R e_3\rightarrow -e_3$.
This in turn implies that $\bar R$ converges to one of the four matrices $R_\alpha R_i^* R_\alpha^T$  $(i=0,\cdots, 3)$, with the first one equal to $I$.

We now prove Property 2. It is straightforward to verify that the last two equilibria $(\bar v, \bar R) = (\frac{2g}{k_1^v}e_3, R_\alpha R_{2,3}^* R_\alpha^T)$ are unstable, since the corresponding equilibria of the subsystem $(\bar v, \bar Re_3)$ are unstable (as a result of Proposition \ref{propo:gammaObs1}). Denoting $\eta \triangleq R_\alpha^T \bar R R_\alpha e_1$, one obtains $(\bar  m_{\cal I} \times \bar R \bar  m_{\cal I})^T e_3 = e_2^T \eta$ and verifies from \eqref{dotRbarZeroDyn} that \vspace{-0.2cm}
\begin{equation} \label{doteta}
\dot \eta = -\bar k_2^r (e_2^T\eta)  \, e_3 \times \eta \vspace{-0.2cm}
\end{equation}
The linearized system of \eqref{doteta} about the ``undesired'' equilibrium $\eta = R_1^* e_1 = -e_1$ satisfies \vspace{-0.2cm}
\[
\dot\eta_2 = \bar k_2^r \eta_2,\quad \dot\eta_1 = \dot \eta_3 = 0 \vspace{-0.2cm}
\]
which clearly indicates that this  ``undesired'' equilibrium is unstable. On the other hand, the linearized system of \eqref{doteta} about the ``desired'' equilibrium $\eta = R_0^* e_1 = e_1$ is given by  \vspace{-0.2cm}
\[
\dot\eta_2 = -\bar k_2^r \eta_2,\quad \dot\eta_1 = \dot \eta_3 = 0  \vspace{-0.2cm}
\]
From here, the local exponential stability of this equilibrium is directly deduced (using the zero-dynamics \eqref{dotRbarZeroDyn}). This along with the local exponential stability of $(\bar v, \bar R e_3) = (0,e_3)$ proved in Proposition \ref{propo:gammaObs1} yields the local exponential stability of the ``desired'' equilibrium $(\bar v,\bar R) = (0,I)$. This concludes the proof.
\end{proof}

Finally, similar results for Observer 2, with proof identical to the one of Theorem \ref{theo:Observer1}, can be directly established.

\vspace{-0.3cm}
\begin{theo} \label{theo:Observer2} (Observer 2)-- Consider System \eqref{sys:systemComplet} and Observer 2 (i.e. \eqref{sys:observerGeneral}+\eqref{innovationObs2}). Assume that $m_{\cal I} \times e_3 \neq 0$. Then, all the properties stated in Theorem \ref{theo:Observer1} hold.
\end{theo}

%\vspace{-0.3cm}
\section{Simulation results}  \label{sec:simulation}\vspace{-0.cm}
Simulations are conducted on a model of a ducted-fan VTOL aerial drone, which was also used in \cite{hua10cep}. Details on the vehicle's model are given in \cite{huaThesis}.
The vehicle is controlled by feedback to track a circular reference trajectory, with the linear velocity given by $\dot x_r =[-15\alpha\sin(\alpha t);15\alpha\cos(\alpha t);0]$ $(m/s)$, with $\alpha=2/\sqrt{15}$. The magnitude of the reference linear acceleration is equal to $4 (m/s^2)$. Due to aerodynamic forces acting on the vehicle, its orientation constantly varies in large proportions.
The normalized earth's magnetic field is taken as $m_{\mathcal I}= [0.434; -0.0091; 0.9008]$.\\
The initial conditions are chosen such that the initial error variables are large and satisfy $\tilde v(0) = [-5;5;-5](m/s)$ and $\bar R(0)  = \text{diag}(-1,1,-1)$.
The following gains are chosen  so that \eqref{eq:gainCond} holds:
$k_1^v = 1.2, k_2^v = 1.2, k_1^r = 0.147, k_2^r = 2.764$,
where the values of $k_1^v, k_2^v, k_1^r$ ensure that the linearized system \eqref{eq:sysLinGoodEqui-1} has a double negative real pole equal to $1.2$. The value of $k_2^r$ is chosen such that one pole of \eqref{doteta} is also equal to $1.2$. In the following, two simulations are reported.

\begin{figure}[!t]\centering%
\psfrag{t (s)}{\scriptsize $t (s)$}%
\psfrag{roll (deg)}{\scriptsize $\text{roll}\, \phi \, \text{(deg)}$}%
\psfrag{pitch (deg)}{\scriptsize $\text{pitch}\, \theta \, \text{(deg)}$}%
\psfrag{yaw (deg)}{\scriptsize $\text{yaw}\, \psi \, \text{(deg)}$}%
\psfrag{Real value}{\scriptsize $\text{Real value}$}%
\psfrag{Observer 1}{\scriptsize $\text{Observer 1}$}%
\psfrag{Observer 2}{\scriptsize $\text{Observer 2}$}%
\includegraphics[width=0.985\linewidth]{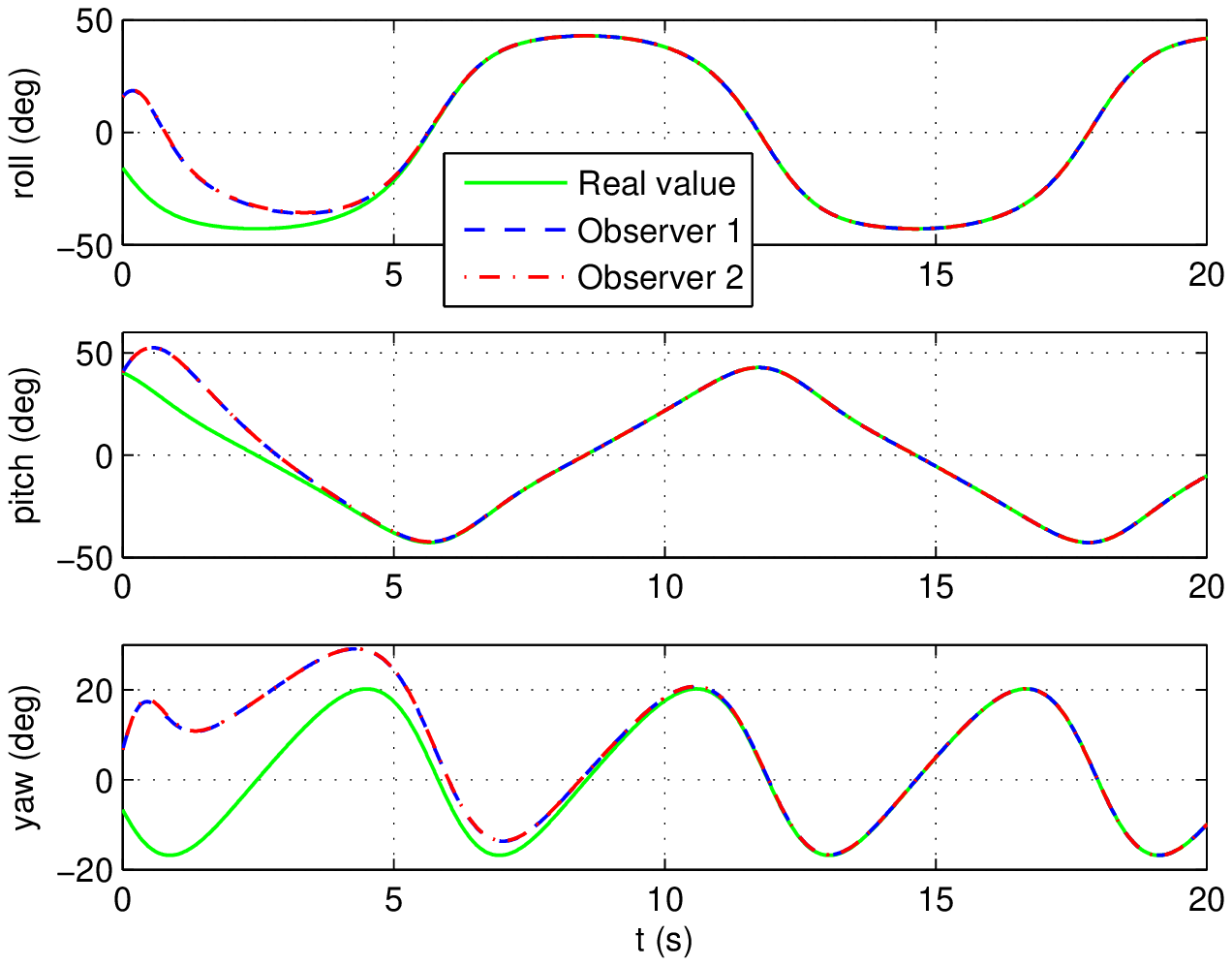} \vspace*{-0.5cm}
\caption{Estimated and real Euler angles versus time (Simulation 1).} \label{fig1Simu1} %\vspace{-0.4cm}
\psfrag{t (s)}{\scriptsize $t (s)$}%
\psfrag{v1 (m/s)}{\scriptsize $v_1 \text{(m/s)}$}%
\psfrag{v2 (m/s)}{\scriptsize $v_2 \text{(m/s)}$}%
\psfrag{v3 (m/s)}{\scriptsize $v_3 \text{(m/s)}$}%
\psfrag{Real value}{\scriptsize $\text{Real value}$}%
\psfrag{Observer 1}{\scriptsize $\text{Observer 1}$}%
\psfrag{Observer 2}{\scriptsize $\text{Observer 2}$}%
\includegraphics[width=0.985\linewidth]{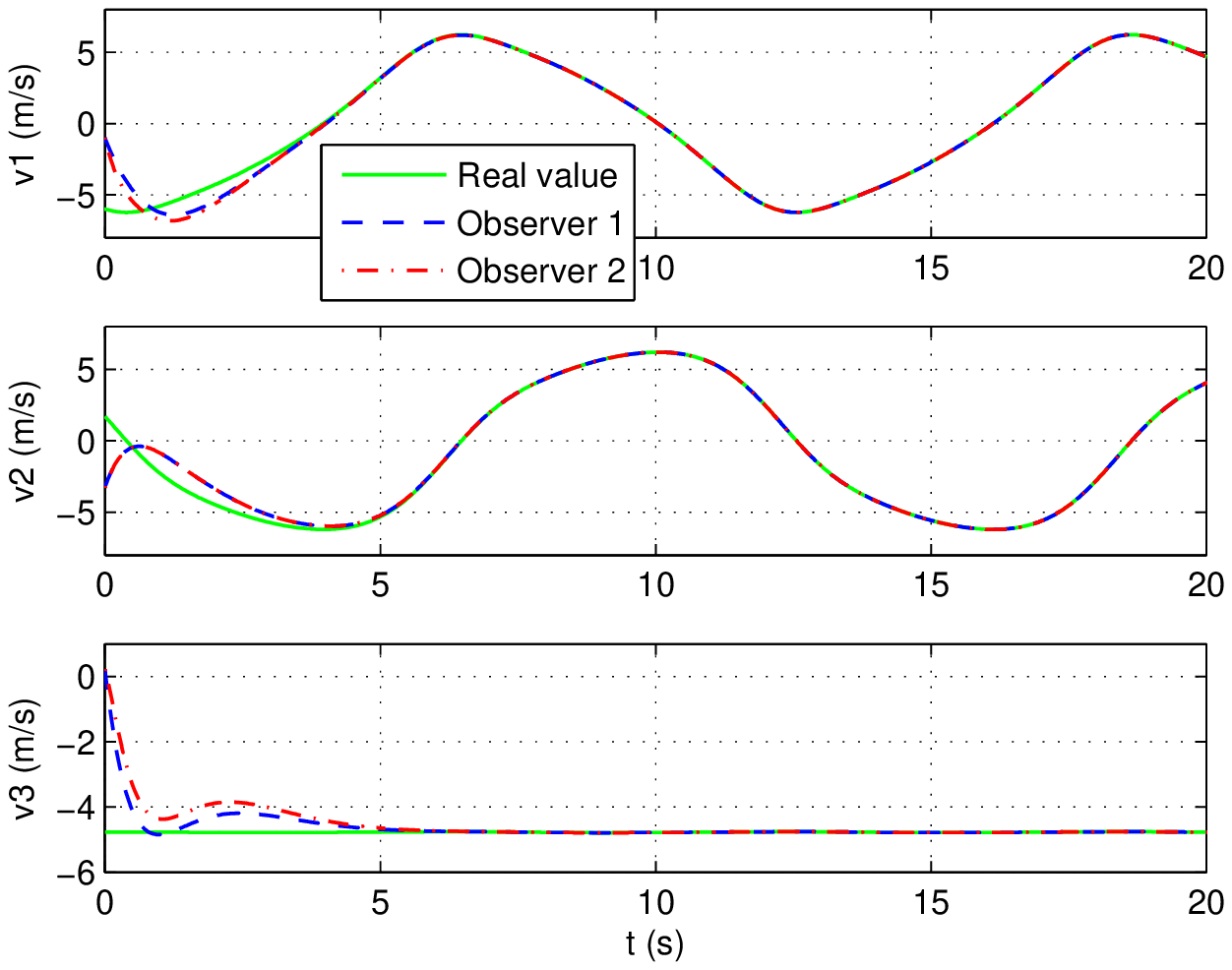} \vspace*{-0.5cm}
\caption{Estimated and real velocity versus time (Simulation 1).} \label{fig2Simu1} \vspace{-0.65cm}
\end{figure}

\noindent $\bullet$  {\em Simulation 1:} This simulation allows us to show the performance of these two observers in the case of perfect measurements. The time evolution of the estimated and real attitudes, represented by Euler angles, along with the estimated and real velocity is shown in Figs.~\ref{fig1Simu1} and \ref{fig2Simu1}, respectively. Both observers ensure the asymptotic convergence of the estimated variables to the real values despite the large initial estimation errors. Their convergence rates are similar and quite satisfactory.

\begin{figure}[!t]\centering%
\psfrag{t (s)}{\scriptsize $t (s)$}%
\psfrag{roll (deg)}{\scriptsize $\text{roll}\, \phi \, \text{(deg)}$}%
\psfrag{pitch (deg)}{\scriptsize $\text{pitch}\, \theta \, \text{(deg)}$}%
\psfrag{yaw (deg)}{\scriptsize $\text{yaw}\, \psi \, \text{(deg)}$}%
\psfrag{Real value}{\scriptsize $\text{Real value}$}%
\psfrag{Observer 1}{\scriptsize $\text{Observer 1}$}%
\psfrag{Observer 2}{\scriptsize $\text{Observer 2}$}%
\includegraphics[width=0.985\linewidth]{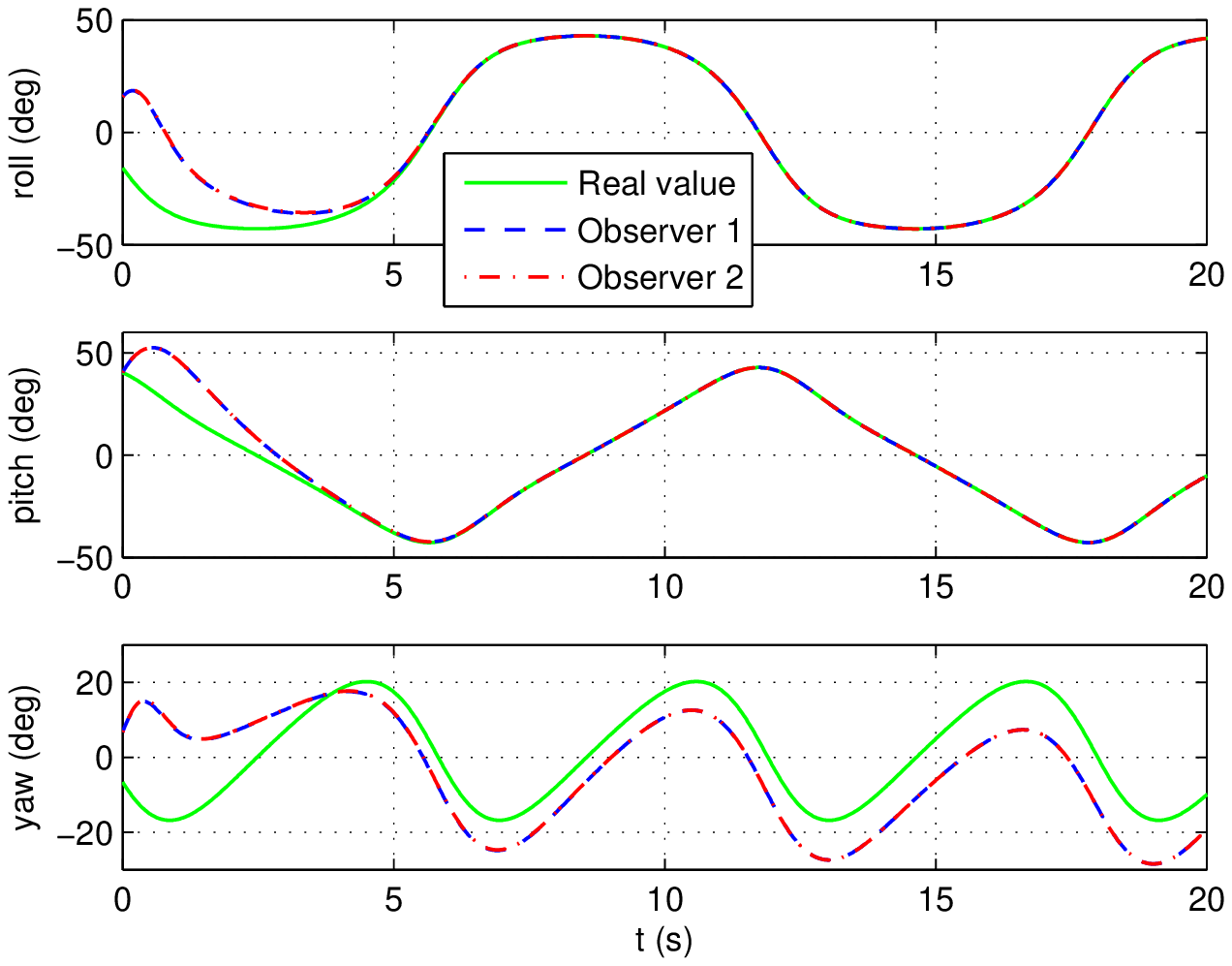} \vspace*{-0.5cm}
\caption{Estimated and real Euler angles versus time (Simulation 2).} \label{fig1Simu2} %\vspace{-0.4cm}
\psfrag{t (s)}{\scriptsize $t (s)$}%
\psfrag{v1 (m/s)}{\scriptsize $v_1 \text{(m/s)}$}%
\psfrag{v2 (m/s)}{\scriptsize $v_2 \text{(m/s)}$}%
\psfrag{v3 (m/s)}{\scriptsize $v_3 \text{(m/s)}$}%
\psfrag{Real value}{\scriptsize $\text{Real value}$}%
\psfrag{Observer 1}{\scriptsize $\text{Observer 1}$}%
\psfrag{Observer 2}{\scriptsize $\text{Observer 2}$}%
\includegraphics[width=0.985\linewidth]{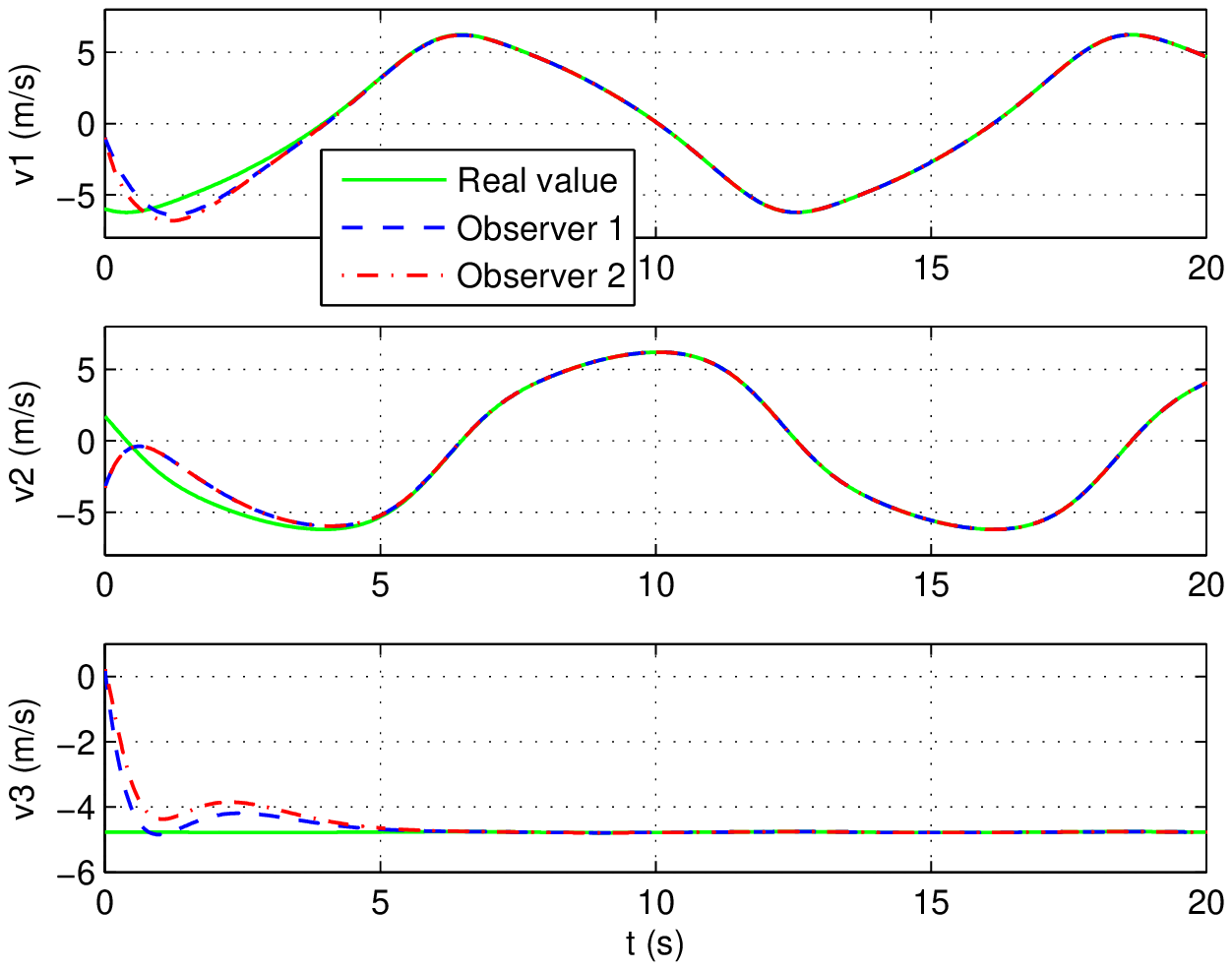} \vspace*{-0.5cm}
\caption{Estimated and real velocity versus time (Simulation 2).} \label{fig2Simu2} \vspace{-0.65cm}
\end{figure}

\noindent $\bullet$  {\em Simulation 2:} With respect to Simulation 1, we only add a constant bias to magnetometer measurements. As expected, it can be observed from Figs.~\ref{fig1Simu2} and \ref{fig2Simu2} that the estimated roll and pitch angles and the estimated velocity components still converge to the real ones, and that the magnetometer measurement bias only affects the yaw estimation (see Fig.~\ref{fig1Simu2}--bottom). This confirms the global decoupling of roll and pitch estimation from magnetometer measurements.

\vspace{-0.1cm}
\section{Conclusion} \label{sec:conclusion}
In this paper, the problem of attitude estimation for accelerated rigid bodies is re-visited, and two novel nonlinear invariant observers are proposed based on the fusion of measurement data provided by an IMU and the measurement of the linear velocity expressed in the body-fixed frame. The paper provides rigourous Lyapunov-based analyses of convergence and stability showing that both observers are almost globally asymptotically stable and locally exponentially stable. Moreover, the roll and pitch estimation is globally decoupled from magnetometer measurements and gain tuning can be easily done, which are interesting from practical standpoints.

\vspace{-0.1cm}
\bibliography{./bibfile}

\begin{thebibliography}{10}

\bibitem{BonnaMR2006ACC}
S.~Bonnabel, P.~Martin, and P.~Rouchon.
\newblock A non-linear symmetry-preserving observer for velocity-aided inertial
  navigation.
\newblock In {\em American Control Conference}, pages 2910--2914, 2006.

\bibitem{BonnabelITAC08}
S.~Bonnabel, P.~Martin, and P.~Rouchon.
\newblock Symmetry-preserving observers.
\newblock {\em IEEE Trans. on Automatic Control}, 53(11):2514--2526, 2008.

\bibitem{bras2013global}
S.~Br{\'a}s, P.~Rosa, C.~Silvestre, and P.~Oliveira.
\newblock Global attitude and gyro bias estimation based on set-valued
  observers.
\newblock {\em Systems \& Control Letters}, 62(10):937--942, 2013.

\bibitem{Crassidis07}
J.~L. Crassidis, F.~L. Markley, and Y.~Cheng.
\newblock Survey of nonlinear attitude estimation methods.
\newblock {\em AIAA Journal of Guidance, Control, and Dynamics}, 30(1):12--28,
  2007.

\bibitem{Dukan13}
F.~Dukan and A.~J. Sorensen.
\newblock {Integration filter for APS, DVL, IMU and pressure gauge for
  underwater vehicles}.
\newblock In {\em IFAC Conference on Control Applications in Marine Systems},
  pages 280--285, 2013.

\bibitem{grip12}
H.~F. Grip, T.~I. Fossen, T.~A. Johansen, and A.~Saberi.
\newblock Attitude estimation using biased gyro and vector measurements with
  time-varying reference vectors.
\newblock {\em IEEE Trans. on Automatic Control}, 57(5):1332–--1338, 2012.

\bibitem{huaThesis}
M.-D. Hua.
\newblock {\em Contributions to the automatic control of aerial vehicles}.
\newblock PhD thesis, Ph. D. thesis, University of Nice-Sophia Antipolis, 2009.

\bibitem{hua10cep}
M.-D. Hua.
\newblock {Attitude estimation for accelerated vehicles using GPS/INS
  measurements}.
\newblock {\em Control Engineering Practice}, 18(7):723--732, 2010.

\bibitem{hua14}
M.-D. Hua, G.~Ducard, T.~Hamel, R.~Mahony, and K.~Rudin.
\newblock Implementation of a nonlinear attitude estimator for aerial robotic
  vehicles.
\newblock {\em IEEE Trans. on Control Systems Technology}, 22(1):201--213,
  2014.

\bibitem{huaIFAC11}
M.-D. Hua, K.~Rudin, G.~Ducard, T.~Hamel, and R.~Mahony.
\newblock Nonlinear attitude estimation with measurement decoupling and
  anti-windup gyro-bias compensation.
\newblock In {\em IFAC World Congress}, pages 2972--2978, 2011.

\bibitem{jensen11}
K.~J. Jensen.
\newblock Generalized nonlinear complementary attitude filter.
\newblock {\em AIAA Journal of Guidance, Control, and Dynamics},
  34(5):1588--–1592, 2011.

\bibitem{mhp08}
R.~Mahony, T.~Hamel, and J.-M. Pflimlin.
\newblock Nonlinear complementary filters on the special orthogonal group.
\newblock {\em IEEE Trans. on Automatic Control}, 53(5):1203--1218, 2008.

\bibitem{ms08IFAC}
P.~Martin and E.~Salaun.
\newblock An invariant observer for {E}arth-{V}elocity-{A}ided attitude heading
  reference systems.
\newblock In {\em IFAC World Congress}, pages 9857--9864, 2008.

\bibitem{ms10cep}
P.~Martin and E.~Salaun.
\newblock {Design and implementation of a low-cost observer-based attitude and
  heading reference system}.
\newblock {\em Control Engineering Practice}, 18(7):712--722, 2010.

\bibitem{micaelli1993}
A.~Micaelli and C.~Samson.
\newblock Trajectory tracking for unicycle-type and two-steering-wheels mobile
  robots.
\newblock Technical Report 2097, INRIA, 1993.

\bibitem{nf99}
H.~Nijmeijer and T.~I. Fossen.
\newblock {\em New directions in nonlinear observer design}, volume 244 of
  Lecture Notes in Control and Information Sciences.
\newblock Springer, 1999.

\bibitem{rt11CDC}
A.~Roberts and A.~Tayebi.
\newblock {On the attitude estimation of accelerating rigid-bodies using GPS
  and IMU measurements}.
\newblock In {\em IEEE Conf. on Decision and Control}, pages 8088–--8093,
  2011.

\bibitem{Salcudean91}
S.~Salcudean.
\newblock A globally convergent angular velocity observer for rigid body
  motion.
\newblock {\em IEEE Trans. on Automatic Control}, 36(12):1493–--1497, 1991.

\bibitem{tmrm07}
A.~Tayebi, S.~McGilvray, A.~Roberts, and M.~Moallem.
\newblock Attitude estimation and stabilization of a rigid body using low-cost
  sensors.
\newblock In {\em IEEE Conf. on Decision and Control}, pages 6424--6429, 2007.

\bibitem{ts03}
J.~Thienel and R.~M. Sanner.
\newblock A coupled nonlinear spacecraft attitude controller and observer with
  an unknown constant gyro bias and gyro noise.
\newblock {\em IEEE Trans. on Automatic Control}, 48(11):2011--2015, 2003.

\bibitem{troniICRA13}
G.~Troni and L.~L. Whitcomb.
\newblock {Preliminary experimental evaluation of a Doppler-aided attitude
  estimator for improved Doppler navigation of underwater vehicles}.
\newblock In {\em IEEE International Conference on Robotics and Automation
  (ICRA)}, pages 4134--4140, 2013.

\bibitem{trumpf12}
J.~Trumpf, R.~Mahony, T.~Hamel, and C.~Lageman.
\newblock Analysis of non-linear attitude observers for time-varying reference
  measurements.
\newblock {\em IEEE Trans. on Automatic Control}, 57:2789--2800, 2012.

\bibitem{vcso08}
J.~F. Vasconcelos, R.~Cunha, C.~Silvestre, and P.~Oliveira.
\newblock A landmark based nonlinear observer for attitude and position
  estimation with bias compensation.
\newblock In {\em IFAC World Congress}, pages 3446–--3451, 2008.

\bibitem{zamani13}
M.~Zamani, J.~Trumpf, and R.~Mahony.
\newblock Minimum-energy filtering for attitude estimation.
\newblock {\em IEEE Trans. on Automatic Control}, 58(11):2917--2921, 2013.

\end{thebibliography}

\end{document}